\documentclass[11pt]{article}
\usepackage{amsfonts, amsmath, amssymb, amsthm}

\PassOptionsToPackage{obeyspaces}{url}
\usepackage[colorlinks=true,citecolor=blue,urlcolor=blue,linkcolor=blue,bookmarksopen=true]{hyperref}
\usepackage{breakurl}
\usepackage{algorithm,algpseudocode,enumitem,tikz}

\author{Boris Bukh\thanks{Department of Mathematical Sciences, Carnegie Mellon University, Pittsburgh, PA 15213. Supported in part by U.S. taxpayers through NSF grant DMS-1301548.} \and Zilin Jiang\footnote{Supported in part by U.S. taxpayers through NSF grant DMS-1201380.}${\ }^{,}$\footnotemark[1]}
\title{A bound on the number of edges in graphs without an even cycle}
\date{}

\newtheorem*{maintheorem}{Main Theorem}

\newtheorem{theorem}{Theorem}
\newtheorem{lemma}[theorem]{Lemma}
\newtheorem{corollary}[theorem]{Corollary}

\newcommand*{\eqdef}{\stackrel{\text{\tiny{def}}}{=}}            
\newcommand*{\abs}[1]{\lvert #1\rvert}                           
\newcommand*{\VBig}{\mathrm{Bg}}
\newcommand*{\VSmall}{\mathrm{Sm}}

\DeclareMathOperator{\ex}{ex}                                    

\def\refwpconds#1{\text{(\hyperref[eq:wpconds]{2#1})}}          

\begin{document}
\maketitle

\begin{abstract}
We show that, for each fixed $k$, an $n$-vertex graph not containing a cycle of length $2k$ has at most
$80\sqrt{k}\log k\cdot n^{1+1/k}+O(n)$ edges.
\end{abstract}

\textsc{MSC classes}: 05C35, 05D99, 05C38

\section*{Introduction}
Let $\ex(n,F)$ be the largest number of edges in an $n$-vertex graph that contains no copy of
a fixed graph~$F$. The systematic study of $\ex(n,F)$ was started by Tur\'an
\cite{turan41} over 70 years ago, and
it has developed into a central problem in extremal graph theory (see surveys
\cite{furedi_simonovits_survey,keevash_survey,sidorenko_survey}).

The function $\ex(n,F)$ exhibits a dichotomy: if $F$ is not bipartite, then $\ex(n,F)$ grows quadratically in $n$,
and is fairly well understood. If $F$ is bipartite, $\ex(n,F)$ is subquadratic, and for very few $F$
the order of magnitude is known. The simplest classes of bipartite graphs are
trees, complete bipartite graphs,
and cycles of even length. Most of the study of $\ex(n,F)$ for bipartite $F$ has been concentrated on these two classes.
In this paper, we address the even cycles. For an overview of the status of $\ex(n,F)$ for complete bipartite graphs see \cite{bbk}.
For a thorough survey on bipartite Tur\'an problems see \cite{furedi_simonovits_survey}.

The first bound on the problem is due to Erd\H{o}s \cite{erdos_integers} who
showed that $\ex(n,C_4)=O(n^{3/2})$. Thanks to the works
of Erd\H{o}s and R\'enyi \cite{erdos_renyi}, Brown \cite[Section~3]{brown}, and K\"ovari, S\'os and Tur\'an \cite{kovari_sos_turan}
it is now known that
\[
  \ex(n,C_4)=(1/2+o(1))n^{3/2}.
\]
The current best bounds for $\ex(n,C_6)$ for large values of $n$ are
\[
  0.5338 n^{4/3}<\ex(n,C_6)\leq 0.6272 n^{4/3}
\]
due to F\"uredi, Naor and Verstra\"ete \cite{fnv}.

A general bound of $\ex(n,C_{2k})\leq \gamma_k n^{1+1/k}$, for some unspecified constant $\gamma_k$, was
asserted by Erd\H{o}s \cite[p.~33]{erdos_cycle}. The first proof was by Bondy and
Simonovits \cite[Lemma 2]{bondy_simonovits}, who
showed that $\ex(n,C_{2k})\leq 20kn^{1+1/k}$ for all sufficiently large $n$. This was improved by
Verstra{\"e}te \cite{verstraete} to $8(k-1)n^{1+1/k}$ and by Pikhurko \cite{pikhurko} to $(k-1)n^{1+1/k}+O(n)$.
The principal result of the present paper is an improvement of these bounds.
\begin{maintheorem}\label{thm:mainthm}
Suppose $G$ is an $n$-vertex graph that contains no $C_{2k}$, and $n\geq (2k)^{8k^2}$ then
\[
  \ex(n, C_{2k})\leq 80\sqrt{k}\log k\cdot n^{1+1/k}+10k^2n.
\]
\end{maintheorem}
In the published version of this paper, the same result was claimed with $\log k$ replaced by $\sqrt{\log k}$.
This is due to a a mistake in verifying condition \refwpconds{e}, which was discovered by Xizhi Liu.

It is our duty to point out that the improvement offered by the Main Theorem is
of uncertain value because we still do not know if $\Theta(n^{1+1/k})$ is the correct
order of magnitude for $\ex(n,C_{2k})$. Only for $k=2,3,5$ are constructions
of $C_{2k}$-free graphs with $\Omega(n^{1+1/k})$ edges known \cite{benson,wenger,lazebnik_ustimenko,mellinger_mubayi}.
The first author believes it to be likely that $\ex(n,C_{2k})=o(n^{1+1/k})$ for all large $k$.
We stress again that the situation is completely different for odd cycles, where
the value of $\ex(n,C_{2k+1})$ is known exactly for all large $n$ \cite{simonovits_method}.

\paragraph{Proof method and organization of the paper} Our proof is
inspired by that of Pikhurko \cite{pikhurko}. Apart from
a couple of lemmas that we quote from \cite{pikhurko}, the present paper
is self-contained. However, we advise the reader to at least
skim \cite{pikhurko} to see the main idea in a simpler setting.

Pikhurko's proof builds a breadth-first search tree, and then argues
that a pair of adjacent levels of the tree cannot contain a $\Theta$-graph\footnote{We recall the definition of a $\Theta$-graph in Section~\ref{sec:theta}}.
It is then deduced that each level
must be at least $\delta/(k-1)$ times larger than the previous, where $\delta$
is the (minimum) degree. The bound on $\ex(n,C_{2k})$ then follows.
The estimate of $\delta/(k-1)$ is sharp when one restricts one's
attention to a pair of levels.

In our proof, we use three adjacent levels.
We find a $\Theta$-graph satisfying an extra technical condition
that  permits an extension of Pikhurko's argument. Annoyingly, this
extension requires a bound on the \emph{maximum degree}.
To achieve such a bound we use a modification of breadth-first
search that avoids the high-degree vertices.

What we really prove in this paper is the following.
\begin{theorem}\label{thm:real}
Suppose $k\geq 4$, and suppose $G$ is a biparite $n$-vertex graph of minimum degree at least $2d+5k^2$,
where
\begin{equation}\label{eq:defd}
  d\geq \max(20\sqrt{k} \log k\cdot n^{1/k},(2k)^{8k}),
\end{equation}
then $G$ contains $C_{2k}$.
\end{theorem}
The Main Theorem follows from Theorem~\ref{thm:real} and two well-known
facts: every graph contains a bipartite subgraph with half of the edges, and
every graph of average degree $d_{\text{avg}}$ contains a subgraph
of minimum degree at least~$d_{\text{avg}}/2$.

The rest of the paper is organized as follows. We present
our modification of breadth-first search in Section~\ref{sec:explore}.
In Section~\ref{sec:theta}, which is the heart of the paper, we explain how to find $\Theta$-graphs
in triples of consecutive levels. Finally, in Section~\ref{sec:final}
we assemble the pieces of the proof.

\section{Graph exploration}\label{sec:explore}
Our aim is to have vertices of degree at most $\Delta d$ for some $k\ll \Delta\ll d^{1/k}$.
The particular choice is fairly flexible; we choose to use
\[
  \Delta\eqdef k^3.
\]

Let $G$ be a graph, and let $x$ be any vertex of $G$. We start our exploration with the set $V_0=\{x\}$, and mark the
vertex $x$ as explored. Suppose $V_0,V_1,\dotsc,V_{i-1}$ are the sets explored in the $0$th, $1$st,\dots,$(i-1)$st steps
respectively. We then define $V_i$ as follows:
\begin{enumerate}
\item Let $V_i'$ consist of those neighbors of $V_{i-1}$ that have not yet been explored.
Let $\VBig_i$ be the set of those vertices in $V_i'$ that have more than $\Delta d$ unexplored neighbors,
and let $\VSmall_i=V_i'\setminus \VBig_i$.
\item Define
\[
  V_i=
   \begin{cases}
     V_i'&\text{if }\abs{\VBig_i}>\tfrac{1}{2k}\abs{V_i'},\\
     \VSmall_i&\text{if } \abs{\VBig_i}\leq \tfrac{1}{2k}\abs{V_i'}.
   \end{cases}
\]
The vertices of $V_i$ are then marked as explored.
\end{enumerate}

We call sets $V_0,V_1,\dotsc$ \emph{levels} of $G$. A level $V_i$ is \emph{big}
if $\abs{\VBig_i}>\tfrac{1}{2k}\abs{V_i'}$, and is \emph{normal} otherwise.

\begin{lemma}\label{lem:mindeg}
If $\delta\leq \Delta d$, and $G$ is a bipartite graph of minimum degree at least $\delta$, then
each $v\in V_{i+1}$ has at least $\delta$ neighbors in $V_i\cup V_{i+2}'$.
\end{lemma}
\begin{proof}
Fix a vertex $v\in V(G)$. We will show, by induction on $i$, that if $v\not\in
V_1\cup\dotsb\cup V_i$, then $v$ has at least $\delta$ neighbors in
$V(G)\setminus (V_1\cup \dotsb\cup V_{i-1})$. The base case $i=1$ is clear.
Suppose $i>1$. If $v\in \VBig_i$, then $v$ has $\Delta d\geq \delta$ neighbors
in the required set. Otherwise, $v$ is not in $V_i'$ and hence has
no neighbors in $V_{i-1}$. Hence, $v$ has as many neighbors in
$V(G)\setminus (V_1\cup \dotsb\cup V_{i-1})$ as in $V(G)\setminus (V_1\cup \dotsb\cup V_{i-2})$,
and our claim follows from the induction hypothesis.

If $v\in V_{i+1}$, then the neighbors of $v$ are a subset of $V_1\cup\dotsb\cup V_i\cup V_{i+2}'$.
Hence, at least $\delta$ of these neighbors lie in $V_i\cup V_{i+2}'$.
\end{proof}

\paragraph{Trilayered graphs}
We abstract out the properties of a triple of consecutive levels into the following definition.
A \emph{trilayered graph} with layers $V_1,V_2,V_3$ is a graph $G$ on a vertex set $V_1\cup V_2\cup V_3$ such that the
only edges in $G$ are between $V_1$ and $V_2$, and between $V_2$ and $V_3$. If $V_1'\subset V_1$, $V_2'\subset V_2$
and $V_3'\subset V_3$, then we denote by $G[V_1',V_2',V_3']$ the trilayered subgraph induced by three sets $V_1',V_2',V_3'$.
Because the graph $G$ that has been explored is bipartite, there are no edges inside each level. Therefore any three sets $V_{i-1},V_i,V_{i+1}'$ from the exploration process naturally form a trilayered graph; these graphs
and their subgraphs are the only trilayered graphs that appear in this paper.

\parshape=13 0cm.82\hsize 0cm.82\hsize 0cm.82\hsize 0cm.82\hsize 0cm.82\hsize 0cm.82\hsize 0cm.82\hsize 0cm.82\hsize 0cm.82\hsize 0cm.82\hsize 0cm.82\hsize 0cm\hsize 0cm \hsize
We say that a trilayered graph has \emph{minimum degree}\vadjust{\hfill\smash{\lower 70pt\llap{%
\begin{tikzpicture}
\draw (0.0,0) ellipse (0.3 and 0.8);
\draw (0.8,0) ellipse (0.3 and 0.8);
\draw (1.6,0) ellipse (0.3 and 0.8);
\node at (0.0, 1) {$V_1$};
\node at (0.8, 1) {$V_2$};
\node at (1.6, 1) {$V_3$};
\node[align=center] at (0.0,-1.2) {$\to$\\[-0.7ex]$A$};
\node[align=center] at (0.59,-1.2) {$\gets$\\[-0.7ex]$B$};
\node[align=center] at (1.01,-1.2) {$\to$\\[-0.7ex]$C$};
\node[align=center] at (1.6,-1.2) {$\gets$\\[-0.7ex]$D$};
\end{tikzpicture}}}}
at least $[A:B,C:D]$ if each vertex in $V_1$ has at least $A$ neighbors in $V_2$, each vertex in $V_2$ has at
least $B$ neighbors in $V_1$, each vertex in $V_2$ has at least $C$ neighbors in $V_3$, and each vertex in $V_3$
has at least $D$ neighbors in $V_2$. A schematic drawing of such a graph is on the right.

\section{\texorpdfstring{$\Theta$}{Theta}-graphs}\label{sec:theta}
A \emph{$\Theta$-graph} is a cycle of length at least $2k$ with a chord.
We shall use several lemmas from the previous works.
\begin{lemma}[Lemma~2.1 in \cite{pikhurko}, also Lemma~2 in \cite{verstraete}]\label{lem:bipartitetheta}
Let $F$ be a $\Theta$-graph and $1\leq l\leq \abs{V(F)}-1$. Let $V(F)=W\cup Z$ be an arbitrary
partition of its vertex set into two non-empty parts such that every path in $F$ of length
$l$ that begins in $W$ necessarily ends in~$W$. Then $F$ is bipartite with parts $W$ and~$Z$.
\end{lemma}
\begin{lemma}[Lemma~2.2 in \cite{pikhurko}]\label{lem:avg_bipartite}
Let $k\geq 3$. Any bipartite graph $H$ of minimum degree at least $k$ contains a $\Theta$-graph.
\end{lemma}
\begin{corollary}\label{cor:avg_bipartite}
Let $k\geq 3$. Any bipartite graph $H$ of average degree at least $2k$ contains a $\Theta$-graph.
\end{corollary}

For a graph $G$ and a set $Y\subset V(G)$ let $G[Y]$ denote the graph induced on $Y$.
For disjoint $Y,Z\subset V(G)$ let $G[Y,Z]$ denote the bipartite subgraph of $G$ that is induced
by the bipartition $Y\cup Z$.

\paragraph{Well-placed $\Theta$-graphs}
Suppose $G$ is a trilayered graph with layers $V_1,V_2,V_3$. We say that a \mbox{$\Theta$-graph} $F\subset G$
is \emph{well-placed} if each vertex of $V(F)\cap V_2$ is adjacent to some vertex in $V_1\setminus V(F)$.
The condition ensures that, for each vertex $v$ of $F$ in $V_2$ there exists a path from the root to $v$
that avoids $F$.

\begin{lemma}\label{lem:trilayered}
Suppose $G$ is a trilayered graph with layers $V_1$, $V_2$, $V_3$ such that
the degree of every
vertex in $V_2$ is at least $2d+5k^2$, and no vertex in $V_2$ has more than $\Delta d$ neighbors in $V_3$.
Suppose $t$ is a nonnegative integer, and let
$F=\frac{d\cdot e(V_1,V_2)}{8k\abs{V_3}}$. Assume that
\begin{equation}
\label{eq:wpconds}
\begin{aligned}
\text{a) }&&F&\geq 2,\\
\text{b) }&&  e(V_1,V_2)&\geq 2k F  \abs{V_1},\\
\text{c) }&&  e(V_1,V_2)&\geq 8k(t+1)^2(2\Delta k)^{2k-1}\abs{V_1},\\
\text{d) }&&  e(V_1,V_2)&\geq 8 (et/F)^t k\abs{V_2},\\
\text{e) }&&  e(V_1,V_2)&\geq 20(t+1)^2\abs{V_2}.
\end{aligned}
\end{equation}

Then at least one of the following holds:
\begin{enumerate}[label=\Roman{*}),ref=(\Roman{*})]
\item There is a $\Theta$-graph in $G[V_1,V_2]$.
\item There is a well-placed $\Theta$-graph in $G[V_1,V_2,V_3]$.
\end{enumerate}
\end{lemma}
The proof of Lemma~\ref{lem:trilayered} is in two parts: finding trilayered subgraph of large
minimum degree (Lemmas \ref{lem:subtrilayer} and \ref{lem:subtrilayeriterated}), and finding
a well-placed $\Theta$-graph inside that trilayered graph (Lemma~\ref{lem:min_trilayered}).

\paragraph{Finding a trilayered subgraph of large minimum degree}
The disjoint union of two bipartite graphs shows that
a trilayered graph with many edges need not contain a trilayered subgraph of large minimum degree.
We show that, in contrast, if a trilayered
graph contains no $\Theta$-graph between two of its levels, then it must contain a subgraph
of large minimum degree. The next lemma demonstrates a weaker version of this claim:
it leaves open a possibility that the graph contains a denser trilayered subgraph.
In that case, we can iterate inside that subgraph, which is done in Lemma~\ref{lem:subtrilayeriterated}.

\begin{lemma}\label{lem:subtrilayer}
Let $a,A,B,C,D$ be positive real numbers. Suppose $G$ is a trilayered graph with layers $V_1$, $V_2$, $V_3$
and the degree of every vertex in $V_2$ is at least $d+4k^2+C$. Assume
also that
\begin{equation}\label{eq:subtrilayercond}
  a\cdot e(V_1,V_2)\geq (A+k+1)\abs{V_1}+B\abs{V_2}.
\end{equation}

Then one of the following holds:

\begin{enumerate}[label=\Roman{*}),ref=(\Roman{*})]
\item \label{alt:v1v2theta} There is a $\Theta$-graph in $G[V_1,V_2]$.

\item \label{alt:trilayer} There exist non-empty subsets $V_1'\subset V_1$, $V_2'\subset V_2$, $V_3'\subset V_3$
such that the induced trilayered subgraph $G[V_1',V_2',V_3']$ has minimum degree at least
$[A:B,C:D]$.

\item \label{alt:increment} There is a subset $\widetilde{V}_2\subset V_2$ such that $e(V_1,\widetilde{V}_2)\geq (1-a)e(V_1,V_2)$,
and $\abs{\widetilde{V}_2}\leq D\abs{V_3}/d$.

\end{enumerate}
\end{lemma}
\begin{proof}
We suppose that alternative \ref{alt:v1v2theta} does not hold. Then,
by Corollary~\ref{cor:avg_bipartite}, the average degree of
every subgraph of $G[V_1,V_2]$ is at most $2k$.

Consider the process that aims to construct a subgraph satisfying \ref{alt:trilayer}.
The process starts with $V_1'=V_1$, $V_2'=V_2$ and $V_3'=V_3$,
and at each step removes one of the vertices that violate the minimum degree condition
on $G[V_1',V_2',V_3']$.
The process stops when either no vertices are left, or the minimum degree
of $G[V_1',V_2',V_3']$ is at least $[A:B,C:D]$. Since in the latter case we are done, we assume that
this process eventually removes every vertex of $G$.

Let $R$ be the vertices of $V_2$ that were removed because at the time of removal
they had fewer than $C$ neighbors in $V_3'$. Put
\begin{align*}
  E'&\eqdef \{ uv\in E(G) : u\in V_2,\ v\in V_3,\ \text{and }v\text{ was removed before  }u\},\\
  S&\eqdef \{v \in V_2 : v\text{ has at least }4k^2\text{ neighbors in }V_1\}.
\end{align*}
Note that $\abs{E'}\leq D\abs{V_3}$. We cannot have $\abs{S}\geq \abs{V_1}/k$, for
otherwise the average degree of the bipartite graph $G[V_1,S]$ would be at least $\frac{4k}{1+1/k}\geq 2k$.
So $\abs{S}\leq \abs{V_1}/k$.

The average degree condition on $G[V_1,S]$ implies that
\[
  e(V_1,S)\leq k(\abs{V_1}+\abs{S})\leq (k+1)\abs{V_1}.
\]

Let $u$ be any vertex in $R\setminus S$. Since it is connected to
at least $(d+4k^2+C)-4k^2 = d+C$ vertices of $V_3$, it must be adjacent to at least $d$
edges of $E'$. Thus,
\[
  \abs{R\setminus S}\leq \abs{E'}/d\leq D\abs{V_3}/d.
\]

Assume that the conclusion \ref{alt:increment} does not hold with $\widetilde{V}_2=R\setminus S$. Then
$e(V_1,R\setminus S)<(1-a)e(V_1,V_2)$. Since the total
number of edges between $V_1$ and $V_2$ that were removed
due to the minimal degree conditions on $V_1$ and $V_2$
is at most
$A\abs{V_1}$ and $B\abs{V_2}$ respectively, we conclude that
\begin{eqnarray*}
  e(V_1,V_2)&\leq& e(V_1,S)+e(V_1,R\setminus S)+ A\abs{V_1}+B\abs{V_2}\\
            &<& (k+1)\abs{V_1}+(1-a)e(V_1,V_2)+A\abs{V_1}+B\abs{V_2},
\end{eqnarray*}
implying that
$$a\cdot e(V_1,V_2) < (A+k+1)\abs{V_1}+B\abs{V_2}.$$
The contradiction with \eqref{eq:subtrilayercond} completes the proof.
\end{proof}

\textit{Remark.} The next lemma can be eliminated at the cost of obtaining the bound
$\ex(n,C_{2k})=O(k^{2/3}n^{1+1/k})$ in place of $\ex(n,C_{2k})=O(\sqrt{k}\log k\cdot n^{1+1/k})$.
To do that, we can set $B\approx k^{2/3}$, $D\approx k^{1/3}$ and $a=1/2$.
One can then show that when applied to trilayered graphs arising from the exploration process
the alternative \ref{alt:increment} leads
to a subgraph of average degree $2k$. The two remaining alternatives are dealt by Corollary~\ref{cor:avg_bipartite}
and Lemma~\ref{lem:min_trilayered}.
However, it is possible to obtain a better bound by iterating the preceding lemma.

\begin{lemma}\label{lem:subtrilayeriterated}
Let $C$ be a positive real number.
Suppose $G$ is a trilayered graph with layers $V_1$, $V_2$, $V_3$, and the degree of every vertex in $V_2$ is at least $d+4k^2+C$.
Let $F=\frac{d\cdot e(V_1,V_2)}{8k\abs{V_3}}$, and assume that $F$ and $e(V_1,V_2)$ satisfy \eqref{eq:wpconds} for some integer $t\geq 1$.
Then one of the following holds:

\begin{enumerate}[label=\Roman{*}),ref=(\Roman{*})]
\item \label{altc:v1v2theta} There is a $\Theta$-graph in $G[V_1,V_2]$.

\item \label{altc:trilayer} There exist numbers $A,B,D$ and non-empty subsets $V_1'\subset V_1$, $V_2'\subset V_2$, $V_3'\subset V_3$
such that the induced trilayered subgraph $G[V_1',V_2',V_3']$ has minimum degree at least
$[A:B,C:D]$, with the following inequalities that bind $A$, $B$, and $D$:
\begin{equation}\label{ineq:ABDcond}
\begin{aligned}
  B&\geq 5,\qquad  (B-4)D\geq 2k,\\
  A&\geq 2k(\Delta D)^{D-1}.
\end{aligned}
\end{equation}
\end{enumerate}

\end{lemma}
\begin{proof}
Assume, for the sake of contradiction, that neither \ref{altc:v1v2theta} nor \ref{altc:trilayer} hold.
With hindsight, set $a_j=\frac{1}{t-j+1}$ for $j=0,\dotsc,t-1$.
We shall define a sequence of sets $V_2=V_2^{(0)}\supseteq V_2^{(1)}\supseteq\dotsb\supseteq V_2^{(t)}$ inductively. We denote by
\[
  d_i\eqdef e(V_1,V_2^{(i)})/\abs{V_2^{(i)}}
\]
the average degree from $V_2^{(i)}$ into $V_1$. The sequence $V_2^{(0)},V_2^{(1)},\dotsc,V_2^{(t)}$
will be constructed so as to satisfy
\begin{align}
\label{iter:large}
  e(V_1,V_2^{(i+1)})&\geq (1-a_i)e(V_1,V_2^{(i)}),\\
\label{iter:density}
  d_{i+1}&\geq d_i\cdot F a_i\prod_{j=0}^{i}(1-a_j).
\end{align}
Note that \eqref{iter:large} and the choice of $a_0,\dotsc,a_i$ imply that
\begin{equation}\label{eq:lowerboundonedgeset}
  e(V_1,V_2^{(i)})\geq \tfrac{1}{t+1}e(V_1,V_2).
\end{equation}

The sequence starts with
$V_2^{(0)}=V_2$. Assume $V_2^{(i)}$ has been defined. We proceed to define $V_2^{(i+1)}$. Put
\begin{align*}
A&=a_i e(V_1,V_2^{(i)})/2\abs{V_1}-k-1,\\
B&= a_id_i/4+5,\\
D&= \min(2k,8k/a_id_i).
\end{align*}
With help of \eqref{eq:lowerboundonedgeset} and \refwpconds{c} it is easy to check that the inequalities \eqref{ineq:ABDcond} hold
for this choice of constants.

In addition,
\begin{align*}
  (A+k+1)\abs{V_1}+B\abs{V_2^{(i)}}&=\tfrac{3}{4}a_i e(V_1,V_2^{(i)})+5\abs{V_2^{(i)}}
 \\&\stackrel{\refwpconds{e}}{\leq} \tfrac{3}{4}a_i e(V_1,V_2^{(i)})+\tfrac{1}{4(t+1)^2}e(V_1,V_2)
 \\&\stackrel{\eqref{eq:lowerboundonedgeset}}{\leq} a_i e(V_1,V_2^{(i)}).
\end{align*}
So, the condition \eqref{eq:subtrilayercond} of Lemma~\ref{lem:subtrilayer} is satisfied for the graph $G[V_1,V_2^{(i)},V_3]$.
By Lemma~\ref{lem:subtrilayer} there is a subset $V_2^{(i+1)}\subset V_2^{(i)}$
satisfying \eqref{iter:large} and
\begin{align*}
  \abs{V_2^{(i+1)}}&\leq D\abs{V_3}/d.
\end{align*}
Next we show that the set $V_2^{(i+1)}$ satisfies inequality \eqref{iter:density}. Indeed, we have
\begin{align*}
  d_{i+1}&=\frac{e(V_1,V_2^{(i+1)})}{\abs{V_2^{(i+1)}}}\geq \frac{(1-a_i)e(V_1,V_2^{(i)})}{D\abs{V_3}/d}\geq(1-a_i)a_id_i\frac{d}{8k\abs{V_3}}e(V_1,V_2^{(i)})\\
        &\stackrel{\eqref{iter:large}}{\geq} (1-a_i)a_id_i\frac{d\cdot e(V_1,V_2)}{8k\abs{V_3}}\prod_{j=0}^{i-1}(1-a_j)
=d_i\cdot F a_i\prod_{j=0}^{i}(1-a_j).
\end{align*}\smallskip

Iterative application of \eqref{iter:density} implies
\begin{equation}\label{ineq:dt}
  d_t\geq d_0 F^t \prod_{j=0}^{t-1} a_j(1-a_j)^{t-j}\geq d_0 F^t \prod_{j=0}^{t-1} \frac{e^{-1}}{t-j+1}=d_0 \frac{(F/e)^t}{(t+1)!}.
\end{equation}\smallskip

If we have $\abs{V_2^{(t)}}<\abs{V_1}$, then the average degree of induced subgraph $G[V_1,V_2^{(t)}]$ is
greater than $e(V_1,V_2^{(t)})/\abs{V_1}\stackrel{\eqref{eq:lowerboundonedgeset}}{\geq} e(V_1,V_2)/(t+1)\abs{V_1}\stackrel{\refwpconds{c}}{\geq} 2k$, which by Corollary~\ref{cor:avg_bipartite} leads to outcome \ref{altc:v1v2theta}.

If $\abs{V_2^{(t)}}\geq \abs{V_1}$ and $d_t\geq 4k$, then the average degree of
$G[V_1,V_2^{(t)}]$ is at least $d_t/2\geq 2k$ because $d_t$ is the average degree of $V_2^{(t)}$ into $V_1$, again leading to the outcome \ref{altc:v1v2theta}.
So, we may assume that $d_t<4k$. Since $(t+1)!\leq 2t^t$ we deduce from
\eqref{ineq:dt} that
\[
  d_0 < 4k(t+1)!(e/F)^t\leq 8k (et/F)^t.
\]
This contradicts \refwpconds{d}, and so the proof is complete.
\end{proof}

\paragraph{Locating well-placed $\Theta$-graphs in trilayered graphs} We come to the central argument
of the paper. It shows how to embed well-placed $\Theta$-graphs
into trilayered graphs of large minimum degree. Or rather, it shows
how to embed well-placed $\Theta$-graphs into regular trilayered graphs;
the contortions of the previous two lemmas, and the factor of $\log k$ in
the final bound, come from authors' inability to deal with irregular graphs.

\begin{lemma}\label{lem:min_trilayered}
Let $A,B,D$ be positive real numbers. Let $G$ be a trilayered graph with layers
$V_1$, $V_2$, $V_3$ of minimum degree at least $[A:B,d+k:D]$.
Suppose that no vertex in $V_2$ has more than $\Delta d$ neighbors in $V_3$.
Assume also that
\begin{gather}
  B\geq 5\\
  (B-4)D\geq 2k-2\label{cond:bd}\\
  A\geq 2k(\Delta D)^{D-1}.
\end{gather}
Then $G$ contains a well-placed $\Theta$-graph.
\end{lemma}
\begin{proof}
Assume, for the sake of contradiction, that $G$ contains no well-placed $\Theta$-graphs.
Leaning on this assumption we shall build an arbitrary long path $P$ of the form
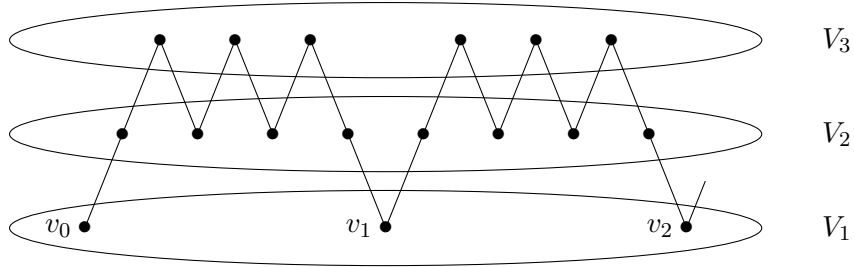
\begin{center}
\begin{tikzpicture}[label distance=-0.5em]
\draw (0,0.00) ellipse (5 and 0.5);
\draw (0,1.25) ellipse (5 and 0.5);
\draw (0,2.50) ellipse (5 and 0.5);
\node at (6,0) {$V_1$};
\node at (6,1.25) {$V_2$};
\node at (6,2.50) {$V_3$};
\draw (-4,0.00) node[label=left:$v_0$] {$\bullet$}
      -- (-3.5,1.25) node {$\bullet$} -- (-3.0,2.50) node {$\bullet$}
      -- (-2.5,1.25) node {$\bullet$} -- (-2.0,2.50) node {$\bullet$}
      -- (-1.5,1.25) node {$\bullet$} -- (-1.0,2.50) node {$\bullet$}
      -- (-0.5,1.25) node {$\bullet$} --
      (0.0,0.0) node[label=left:$v_1$] {$\bullet$}
      -- (+0.5,1.25) node {$\bullet$} -- (+1.0,2.50) node {$\bullet$}
      -- (+1.5,1.25) node {$\bullet$} -- (+2.0,2.50) node {$\bullet$}
      -- (+2.5,1.25) node {$\bullet$} -- (+3.0,2.50) node {$\bullet$}
      -- (+3.5,1.25) node {$\bullet$} --
      (4.0,0.0) node[label=left:$v_2$] {$\bullet$} -- (4.25,0.625);
\end{tikzpicture}
\end{center}
where, for each $i$, vertices $v_i$ and $v_{i+1}$ are joined by a path of length $2D$
that alternates between $V_2$ and $V_3$. Since the graph is finite, this would be a
contradiction.

While building the path we maintain the following property:
\begin{equation}
\label{pc:v2}\tag{$\star$}\text{Every $v\in P\cap V_2$ has at least one neighbor in $V_1\setminus P$}.
\end{equation}
We call a path satisfying \eqref{pc:v2} \emph{good}.

We construct the path inductively. We begin by picking $v_0$ arbitrarily from $V_1$. Suppose a good path $P=v_0\leftrightsquigarrow v_1\leftrightsquigarrow \dotsb \leftrightsquigarrow v_{l-1}$
has been constructed, and we wish to find a path extension $v_0\leftrightsquigarrow v_1\leftrightsquigarrow \dotsb \leftrightsquigarrow v_{l-1}\leftrightsquigarrow v_l$.

There are at least about $A$ ways to extend the path by a single vertex. The idea of the following argument shows that many of these extensions can be extended to another vertex, and then another, and so on.

For each $i=1,2,\dotsc,2D-1$ we shall define a family $\mathcal{Q}_i$ of good paths that satisfy
\begin{enumerate}
\item Each path in $\mathcal{Q}_i$ is of the form $v_0\leftrightsquigarrow v_1\leftrightsquigarrow \dotsb \leftrightsquigarrow v_{l-1} \leftrightsquigarrow u$,
where $v_{l-1} \leftrightsquigarrow u$ is a path of length $i$ that alternates between $V_2$ and $V_3$. The vertex
$u$ is called a \emph{terminal} of the path. The set of terminals of the paths in $\mathcal{Q}_i$ is denoted by $T(\mathcal{Q}_i)$. Note that $T(\mathcal{Q}_i)\subset V_2$ for odd $i$ and $T(\mathcal{Q}_i)\subset V_3$ for even $i$.
\item For each $i$, the paths in $\mathcal{Q}_i$ have distinct terminals.
\item For odd-numbered indices, we have the inequality
\begin{equation}\label{eqn:doublestar}
  \abs{\mathcal{Q}_{2i+1}}\geq -3k+ A\left(\frac{1}{\Delta}\right)^i\prod_{j\leq i} \left(1-\frac{j}{D}\right).
\end{equation}
\item For even-numbered indices, we have the inequality
\begin{equation}\label{eqn:star}
  e(T(\mathcal{Q}_{2i}),V_2)\geq d\abs{\mathcal{Q}_{2i-1}}.
\end{equation}
\end{enumerate}

Let
\[
  t\eqdef \left\lceil B/2\right\rceil.
\]\smallskip

We will repeatedly use the following straightforward fact, which we call the \emph{small-degree argument}:
whenever $Q$ is a good path and $u\in V_2\setminus Q$ is adjacent to the terminal of $Q$, then $u$ is
adjacent to fewer than $t$ vertices in $V_1\cap Q$.  Indeed, if vertex $u$ were adjacent to $v_{j_1},v_{j_2},\dotsc,v_{j_t}\in V_1\cap Q$ with $j_1 < j_2 < \ldots < j_k$,
then $v_{j_2} \leftrightsquigarrow u$ (along path $Q$) and the edge
$uv_{j_2}$ would form a cycle of total length at least\linebreak
$
  2D(t-2)+2\geq 2D(B/2-2)+2\stackrel{\eqref{cond:bd}}{\geq}2k
$.
As $uv_{j_3}$ is a chord of the cycle, and $u$ is adjacent to $v_{j_1}$ that is not on the cycle,
that would contradict the assumption that $G$ contains no well-placed $\Theta$-graph.

The set $\mathcal{Q}_1$ consists of all paths of the form $Pu$ for $u\in V_2\setminus P$. Let us check
that the preceding conditions hold for $\mathcal{Q}_1$.
Vertex $v_{l-1}$ cannot be adjacent to $k$ or more vertices in $P\cap V_2$, for otherwise $G$
would contain a well-placed $\Theta$-graph with a chord through $v_{l-1}$. So, $\abs{\mathcal{Q}_1}\geq A-k$.
Next, consider any $u\in V_2\setminus P$
that is a neighbor of $v_{l-1}$. By the small-degree argument vertex $u$ cannot be adjacent to $t$
or more vertices of $P\cap V_1$, and $Pu$ is good.
\medskip

Suppose $\mathcal{Q}_{2i-1}$ has been defined, and we wish to define $\mathcal{Q}_{2i}$. Consider an arbitrary
path\linebreak $Q=v_0\leftrightsquigarrow v_1\leftrightsquigarrow \dotsb \leftrightsquigarrow v_{l-1} \leftrightsquigarrow u\in \mathcal{Q}_{2i-1}$.
Vertex $u$ cannot have $k$ or more neighbors in $Q\cap V_3$, for otherwise $G$ would contain a well-placed $\Theta$-graph with a chord through $u$.
Hence, there are at least $d$ edges of the form $uw$, where $w\in V_3\setminus Q$. As we vary $u$ we obtain
a family of at least $d\abs{\mathcal{Q}_{2i-1}}$ paths.
We let $\mathcal{Q}_{2i}$ consist of any maximal subfamily of such paths with distinct terminals. The condition \eqref{eqn:star} follows automatically as each vertex of $T(\mathcal{Q}_{2i-1})$ has at least $d$ neighbors in $T(\mathcal{Q}_{2i})$.\medskip

Suppose $\mathcal{Q}_{2i}$ has been defined, and we wish to define $\mathcal{Q}_{2i+1}$. Consider an arbitrary
path\linebreak $Q=v_0\leftrightsquigarrow v_1\leftrightsquigarrow \dotsb \leftrightsquigarrow v_{l-1} \leftrightsquigarrow u\in \mathcal{Q}_{2i}$.
An edge $uw$ is called \emph{long} if $w\in P$, and $w$ is at a distance exceeding $2k$ from $u$ along path~$Q$.
If $uw$ is a long edge, then from $u$ to $Q$ there is only one edge, namely the edge to the predecessor of $u$ on $Q$,
 for otherwise there is a well-placed $\Theta$-graph.
Also, at most $i$ neighbors of $u$ lie on the path $v_{l-1}\leftrightsquigarrow u$. Since $\deg u\geq D$,
it follows that there are at least $(1-i/D)\deg u$ short edges from $u$ that miss $v_{l-1}\leftrightsquigarrow u$.
Thus there is a set $\mathcal{W}$ of at least $(1-i/D)e(T(\mathcal{Q}_{2i}),V_2)$ walks (not necessarily paths!) of the form
$v_0\leftrightsquigarrow v_1\leftrightsquigarrow \dotsb \leftrightsquigarrow v_{l-1} \leftrightsquigarrow uw$
such that $v_{l-1} \leftrightsquigarrow uw$ is a path and $w$ occurs only among the last $2k$
vertices of the walk.

From the maximum degree condition on $V_2$ it follows that walks in $\mathcal{W}$ have at least\linebreak $(1-i/D)e(T(\mathcal{Q}_{2i}),V_2)/\Delta d$ distinct terminals.
A walk fails to be a path only if the terminal vertex lies on $P$. However, since the edge $uw$ is short,
this can happen for at most $2k$ possible terminals. Hence, there is a $\mathcal{Q}_{2i+1}\subset \mathcal{W}$
of size
\begin{equation}\label{eq:qodd}
\abs{\mathcal{Q}_{2i+1}}\geq (1-i/D)e(T(\mathcal{Q}_{2i}),V_2)/\Delta d-2k
\end{equation}
that consists of paths
with distinct terminals. It remains to check that every path in $\mathcal{Q}_{2i+1}$ is good. The only way that
$Q=v_0\leftrightsquigarrow \dotsb \leftrightsquigarrow v_{l-1} \leftrightsquigarrow uw\in \mathcal{Q}_{2i+1}$
may fail to be good is if $w$ has no neighbors in $V_1\setminus Q$. By the small-degree argument
$w$ has fewer than $t$ neighbors in $V_1$. Since $w$ has at least $B$
neighbors in $V_1$ and $B\geq t+2$, we conclude that $w$ has at least \emph{two}
neighbors in $V_1$ outside the path. Of course, the same is true for every terminal
of a path in~$\mathcal{Q}_{2i+1}$. The condition \eqref{eqn:doublestar} for $\mathcal{Q}_{2i+1}$ follows from \eqref{eq:qodd}, \eqref{eqn:star} and from validity of \eqref{eqn:doublestar} for $\mathcal{Q}_{2i-1}$.

\medskip

Note that $\mathcal{Q}_{2D-1}$ is non-empty. Let $Q=v_0\leftrightsquigarrow \dotsb \leftrightsquigarrow v_{l-1} \leftrightsquigarrow u\in \mathcal{Q}_{2D-1}$
be an arbitrary path. Note that since $2D-1$ is odd, $u\in V_2$. By the property of terminals of $V_i$ (odd $i$) that
we noted in the previous paragraph, there are two vertices in $V_1\setminus Q$ that are neighbors of~$u$.
Let $v_l$ be any of them, and let the new path be $Qv_l=v_0\leftrightsquigarrow \dotsb \leftrightsquigarrow v_{l-1} \leftrightsquigarrow uv_l$.
This path can fail to be good if there is a vertex $w$ on the path $Q$ that is good in $Q$, but
is bad in $Qv_l$. By the small-degree argument, $w$ is adjacent to fewer than $t$
vertices in $Q\cap V_1$ that precede $w$ in~$Q$. The same argument applied to the reversal of the path $Qv_l$
shows that $w$ is adjacent to fewer than $t$ vertices in $Q\cap V_1$ that succeeds $w$ in~$Q$.
Since $2t-2<B$, the path $Qv_l$ is good.

Hence, it is possible to build an arbitrarily long path in $G$. This contradicts the finiteness of~$G$.
\end{proof}
Lemma~\ref{lem:trilayered} follows from Lemmas~\ref{lem:subtrilayeriterated} and \ref{lem:min_trilayered}
by setting $C=d+k$, in view of inequality $4k^2+k\leq 5k^2$. We lose $k^2-k$ here for cosmetic reason: $5k^2$ is tidier than $4k^2+k$.

\section{Proof of Theorem~\ref{thm:real}}\label{sec:final}
Suppose that $G$ is a bipartite graph of minimum degree at least $2d+5k^2$ and contains no $C_{2k}$. Pick a root vertex $x$ arbitrarily, and
let $V_0,V_1,\dotsc,V_{k-1}$ be the levels obtained from the exploration process in Section~\ref{sec:explore}.

\begin{lemma}\label{lem:trilayer_notheta}
For $1\leq i\leq k-1$, the graph $G[V_{i-1},V_i,V_{i+1}]$ contains no well-placed $\Theta$-graph.
\end{lemma}
\begin{proof}
The following proof is almost an exact repetition of the proof of Claim~3.1 from \cite{pikhurko} (which is also reproduced
as Lemma~\ref{lem:pikhurko_notheta} below).

Suppose, for the sake of contradiction, that a well-placed $\Theta$-graph $F\subset G[V_{i-1},V_i,V_{i+1}]$ exists.
Let $Y=V_i\cap V(F)$. Since $F$ is well-placed, for every vertex of $Y$ there is a path avoiding $V(F)$ of length $i$
to the vertex~$x$. The union of these paths forms a tree $T$ with $x$ as a root. Let $y$ be the vertex
farthest from $x$ such that every vertex of $Y$ is a $T$-descendant of~$y$. Paths
that connect $x$ to $Y$ branch at $y$. Pick one such branch, and let $W\subset Y$ be the
set of all the $T$-descendants of that branch. Let $Z=V(F)\setminus W$.
From $W\neq V_i\cap V(F)$ it follows that $Z$ is not an independent set of $F$, and so
$W\cup Z$ is not a bipartition of~$F$.

Let $\ell$ be the distance between $x$ and~$y$. We have $\ell<i$ and $2k-2i+2\ell<2k\leq \abs{V(F)}$.
By Lemma~\ref{lem:bipartitetheta} in $F$ there is a path $P$ of length $2k-2i+2\ell$
that starts at some $w\in W$ and ends in $z\in Z$. Since the length of $P$ is even, $z\in Y$.
Let $P_w$ and $P_z$ be unique paths in $T$ that connect $y$ to respectively $w$ and~$z$.
They intersect only at $y$. Each of $P_w$ and $P_z$ has length $i-\ell$. The union of paths
$P,P_w,P_z$ forms a $2k$-cycle in $G$.
\end{proof}

The same argument (with a different $Y$) also proves the next lemma.
\begin{lemma}[Claim~3.1 in \cite{pikhurko}]\label{lem:pikhurko_notheta}
For $1\leq i\leq k-1$, neither of $G[V_i]$ and $G[V_i,V_{i+1}]$ contains a bipartite $\Theta$-graph.
\end{lemma}

The next step is to show that the levels $V_0,V_1,V_2,\dotsc$ increase in size. We shall show by induction on $i$ that
\begin{align}
\label{induc:rightdeg}
  e(V_i,V_{i+1})&\geq d\abs{V_i},\\
\label{induc:leftdeg}
  e(V_i,V_{i+1})&\leq 2 k \abs{V_{i+1}},\\
\label{induc:leftdegprime}
  e(V_i,V_{i+1}')&\leq 2 k \abs{V_{i+1}'},\\
\label{induc:growth}
  \abs{V_{i+1}}&\geq (2k)^{-1} d\abs{V_i},\\
\label{induc:normal}
  \abs{V_{i+1}}&\geq \tfrac{d^2}{400k\log^2 k}\abs{V_{i-1}}.
\end{align}
To prove Theorem~\ref{thm:real}, we only need \eqref{induc:normal}; the remaining inequalities play auxiliary roles in derivation of \eqref{induc:normal}.

Clearly, these inequalities hold for $i=0$ since each vertex of $V_1$ sends only one edge to $V_0$.

\paragraph{Proof of \eqref{induc:rightdeg}:}
By Lemma~\ref{lem:mindeg} the degree of every vertex in $V_i$ is at least
$2d+4k$, and so
\[
  e(V_i,V_{i+1}')\geq (2d+4k)\abs{V_i}-e(V_{i-1},V_i)\stackrel{\text{induc.}}{\geq} (2d+2k)\abs{V_i}.
\]
We next distinguish two cases depending on whether $V_{i+1}$ is big (in the sense of the definition from Section~\ref{sec:explore}).
If $V_{i+1}$ is big, then $e(V_i,V_{i+1})=e(V_i,V_{i+1}')$, and \eqref{induc:rightdeg} follows. If
$V_{i+1}$ is normal, then Corollary \ref{cor:avg_bipartite} and Lemma~\ref{lem:pikhurko_notheta} imply that
\[
  e(V_i,\VBig_{i+1})\leq k(\abs{V_i}+\abs{\VBig_{i+1}})\leq
k\left(\abs{V_i}+\tfrac{1}{2k}\abs{V_{i+1}'}\right) \leq k\abs{V_i} + \tfrac{1}{2}e(V_i,V_{i+1}') \]
and so
\begin{align*}
  e(V_i,V_{i+1})=e(V_i,V_{i+1}')-e(V_i,\VBig_{i+1})\geq \tfrac{1}{2}e(V_i,V_{i+1}')-k\abs{V_i} \geq d\abs{V_i}
\end{align*}
implying \eqref{induc:rightdeg}.\qed

\paragraph{Proof of \eqref{induc:leftdeg}:}
Consider the graph $G[V_i,V_{i+1}]$.
Inequality
\eqref{induc:rightdeg} asserts that the average degree of $V_i$ is at least $d\geq 2k$.
If \eqref{induc:leftdeg} does not hold, then the average degree of $V_{i+1}$ is at least $2k$ as
well, contradicting Corollary~\ref{cor:avg_bipartite} and Lemma~\ref{lem:pikhurko_notheta}.\qed

\paragraph{Proof for \eqref{induc:leftdegprime}:} The argument is the same as for \eqref{induc:leftdeg} with $G[V_i,V_{i+1}']$
in place of $G[V_i,V_{i+1}]$.\qed

\paragraph{Proof for \eqref{induc:growth}:} This follows from \eqref{induc:leftdeg} and \eqref{induc:rightdeg}.\qed

\paragraph{Proof of \eqref{induc:normal} in the case $V_i$ is a normal level:}
We assume that
\eqref{induc:normal} does not hold and will derive a contradiction.
Consider the trilayered graph $G[V_{i-1},V_i,V_{i+1}']$.
Let $t=2\log k$. Suppose momentarily that the inequalities \eqref{eq:wpconds} in Lemma~\ref{lem:trilayered} hold. Then since $V_i$ is normal,
each vertex in $V_i$ has at most $\Delta d$ neighbors in $V_{i+1}'$, and so Lemma~\ref{lem:trilayered} applies. However, the lemma's conclusion contradicts
Lemmas~\ref{lem:trilayer_notheta} and \ref{lem:pikhurko_notheta}. Hence, to prove \eqref{induc:normal} it suffices to
verify inequalities \refwpconds{a--d} with $F=d\cdot e(V_{i-1},V_i)/8k\abs{V_{i+1}'}$.

We may assume that
\begin{equation}\label{eq:Fbound}
F\geq 2e^2\log k,
\end{equation}
and in particular that \refwpconds{a} holds. Indeed, if \eqref{eq:Fbound} were not true, then inequality \eqref{induc:rightdeg} would imply $\abs{V_{i+1}'}\geq (d^2/16e^2k\log k)\abs{V_{i-1}}$,
and thus
\[\abs{V_{i+1}}\geq (1-\tfrac{1}{k})\abs{V_{i+1}'}\geq (d^2/32e^2k\log k)\abs{V_{i-1}},
\]
and so \eqref{induc:normal} would follow in view of $32e^2\leq 400$.

Inequality \refwpconds{b} is implied by \eqref{induc:growth}. Indeed,
\[
  e(V_{i-1},V_i)=8k\abs{V_{i+1}'} F/d \geq 8k\abs{V_{i+1}}F/d \stackrel{\eqref{induc:growth}}{\geq} 4F\abs{V_i}\stackrel{\eqref{induc:growth}}{\geq} 2k^{-1}d F\abs{V_{i-1}},
\]
and $d\geq k^2$ by the definition of $d$ from \eqref{eq:defd}.

Inequality \refwpconds{c} is implied by \eqref{eq:defd} and \eqref{induc:rightdeg}.

Next, suppose \refwpconds{d} were not true. Since $F/t\geq e^2$ by \eqref{eq:Fbound}, we would then conclude
\begin{align*}
\abs{V_{i+1}}&\stackrel{\eqref{induc:growth}}{\geq} (2k)^{-1}d \abs{V_i}{\geq} d(16k^2)^{-1}(F/et)^t e(V_{i-1},V_i)\\&\geq d(16k^2)^{-1} e^{2\log k} e(V_{i-1},V_i)\stackrel{\eqref{induc:rightdeg}}{\geq} \tfrac{1}{16} d^2\abs{V_{i-1}},
\end{align*}
and so \eqref{induc:normal} would follow.

Finally, \refwpconds{e} is a consequence of \eqref{induc:rightdeg}. Indeed, if \refwpconds{e} fails, then
\[
  e(V_{i-1},V_i)\leq 20(2\log k+1)^2\abs{V_i}\stackrel{\eqref{induc:growth}}{\leq} 20(2\log k+1)^2 \frac{2k}{d}\abs{V_{i+1}}\leq 360 \frac{k\log^2 k}{d}\abs{V_{i+1}}.
\]
This inequality and \eqref{induc:rightdeg} then together imply \eqref{induc:normal}.
\qed


\paragraph{Proof of \eqref{induc:normal} in the case $V_i$ is a big level:}
We have
\begin{align*}
  \abs{V_{i+1}}&\geq \tfrac{1}{2}\abs{V_{i+1}'}\stackrel{\eqref{induc:leftdegprime}}{\geq} (4k)^{-1} e(V_i,V_{i+1}')\geq (4k)^{-1}e(\VBig_i,V_{i+1}')\geq (4k)^{-1} \Delta d\abs{\VBig_i}\\
               &\geq (8k^2)^{-1}\Delta d \abs{V_i}\stackrel{\eqref{induc:growth}}{\geq} (16k^3)^{-1}\Delta d^2 \abs{V_{i-1}}=\tfrac{1}{16}d^2\abs{V_{i-1}},
\end{align*}
and so \eqref{induc:normal} holds.\qed

We are ready to complete the proof of Theorem~\ref{thm:real}. If $k$ is even, then  $k/2$ applications of \eqref{induc:normal}
yield
\[
  \abs{V_k}\geq \frac{d^k}{(400k\log^2 k)^{k/2}}.
\]
If $k$ is odd, then $(k-1)/2$ applications of \eqref{induc:normal} yield
\[
  \abs{V_k}\geq \frac{d^{k-1}}{(400k\log^2 k)^{(k-1)/2}}\abs{V_1}\geq \frac{d^k}{(400k\log^2 k)^{(k-1)/2}}.
\]
Either way, since $\abs{V_k}<n$ we conclude that $d<20\sqrt{k}\log k\cdot n^{1/k}$.

\section{Acknowledgment}
We would like to thank the referees for carefully reading the manuscript and for giving constructive comments which helped improving the quality of the paper. We thank  Xizhi Liu for bringing to our attention a mistake in the proof of \eqref{induc:normal} in the original version of the paper, which resulted in us claiming a stronger result with $80\sqrt{k\log k}$ instead of $80\sqrt{k}\log k$. All remaining errors are ours.

\bibliographystyle{plain}
\bibliography{ect}

\appendix
\section{Addendum (joint with Jie Ma)}
After the paper was written and published, we made two observations:
\begin{itemize}
\item The method in the paper cannot improve $80\sqrt{k} \log k$ to anything
better than $O(\sqrt{k})$.
\item In the proof of our main theorem, there is a way to reduce to the case when
$G$ is almost-regular. This will simplify the argument, and might lead to reducing
the power of $\log k$ in the result.
\end{itemize}

\paragraph{Limit of the method:} A fundamental problem in extremal combinatorics is the \emph{girth problem}: to estimate $\ex(n,\{C_3,C_4,\dotsc,C_{2k}\})$, i.e., the size
of the largest graph of girth at least $2k+1$. It is easy to prove that $\ex(n,\{C_3,C_4,\dotsc,C_{2k}\})\leq Cn^{1+1/k}$ for an absolute constant $C$.
Indeed, suppose $G$ is a given graph of girth at least $2k+1$. We pass to a subgraph of a large minimum degree, pick one of the remaining vertices $v$ and consider a depth-first search tree
based at $v$. As all vertices at depth $k$ are distinct, the bound follows\footnote{It is possible to replace minimum degree by average
degree in this sketch. See \cite{alon_hoory_linial}}. All the upper bounds on $\ex(n,C_{2k})$, including ours, are embellishments of this basic argument, as
no other argument for the girth problem is known.

The girth problem admits a generalization to bipartite graphs. Let $\ex(n,m,C_{\leq 2k})$ be the largest number of edges in a
bipartite graph with parts of size $m$ and $n$ of girth at least $2k+1$. The basic argument above easily extends to show
that $\ex(n,m,C_{\leq 2k})\leq Cn^{1/k}\cdot(mn)^{1/2}$ if $k$ is even (and similar, but more complicated expression for odd $k$).
Suppose $k$ is even, and $G$ is a bipartite graph with parts of sizes $n/k$ and $n$ that has $Ck^{-1/2}n^{1+1/k}$ edges. By
cloning each vertex in the smaller part into $k$ copies, we obtain a $C_{2k}$-free $2n$-vertex graph with $Ck^{1/2}n^{1+1/k}$ edges.
So, proving a bound of the from $\ex(n,\{C_3,\dotsc,C_{2k}))=o(\sqrt{k} n^{1+1/k})$ would require improving on the basic girth argument.

A similar construction appears in \cite{naor_verstraete}.

\paragraph{Potential improvement:} Some of technical difficulties in the paper come from
dealing with irregular graphs. It is possible to circumvent them by passing to an almost regular subgraph.
A result of Erd\H{o}s and Simonovits \cite{MR0300924} shows that every sufficiently large $n$-vertex graph with $n^{1+1/k}$ edges contains a subgraph $H$ on $m\geq n^{(1-1/k)/(k+1)}$ vertices with at least
$\tfrac{2}{5}m^{1+1/k}$ edges satisfying
$$\text{maximum degree of }H \le 10\cdot 2^{k^2+1} \cdot \text{minimum degree of }H.$$ It is easy to modify their argument
to handle a graph with $cn^{1+1/k}$ edges instead of $n^{1+1/k}$. A similar result has appeared before in \cite[Proposition~2.7]{jiang_seiver}.

\begin{theorem}
For every $c > 0$ and $\alpha \in (0,1]$, if $n$ is sufficiently large (relative to $c$ and $\alpha$), then every $n$-vertex graph $G$ with $\ge c n^{1+\alpha}$ edges contains a subgraph $G'$ on at least $cn^{\alpha/2}$ vertices such that the degree of each vertex in $G'$ is between $\frac{\alpha}{6}cv(G')^{\alpha}$ and $\frac{2}{\alpha\gamma}cv(G')^{\alpha}$, where $\gamma=6^{-2/\alpha}$.
\end{theorem}

\begin{proof}
There are two parts to the argument. We first find a subgraph in which the ratio between the minimum and maximum degrees is bounded. We then
thin that subgraph on $m$ vertices to make the degrees approximately $cm^{\alpha}$.

\medskip\noindent\textit{Finding dense subgraphs:} Let $H$ be a subgraph of $G$ that maximizes the ratio $e(H) / v(H)^{1+\alpha/2}$. By the assumption on $e(G)$, this ratio is at least $cn^{\alpha/2}$.
Since $e(H)\leq v(H)^2/2$ and $v(H) \le n$, it then follows that $$v(H)^{1-\alpha/2}\geq 2cn^{\alpha/2} \text{ and }e(H) \ge cv(H)^{1+\alpha}.$$
Let $S$ be the subset of $V(H)$ consisting of $\gamma v(H)$ vertices with the largest degrees in $H$.

Suppose first that at least $e(H)/2$ edges of $H$ are incident to a vertex in $S$. Set $\eta=\gamma/\alpha$. By averaging, we can find a set $T\subset V(H)\setminus S$ of $\eta v(H)$ vertices that is incident to at least $\eta/(1-\gamma)$ fraction of edges of $H$ from $S$ to~$V(H)\setminus S$. Hence,
$$
e(H[S\cup T]) \geq \left(\frac{\eta}{1-\gamma}\right)e(S,T)+e(S)
              \geq \left(\frac{\eta}{1-\gamma}\right)\bigl(e(S,T)+e(S)\bigr)
              \geq \left(\frac{\eta}{1-\gamma}\right)\frac{e(H)}{2}
                >  \frac{\eta}{2}e(H).
$$
Since \[
(\gamma+\eta)^{1+\alpha/2}=\gamma^{1+\alpha/2}(1+1/\alpha)^{1+\alpha/2}\leq (3/\alpha)\gamma^{1+\alpha/2} = \eta/2,
\]
it follows that $e(H[S\cup T])/v(H[S\cup T])^{1+\alpha/2} > e(H)/v(H)^{1+\alpha/2}$, contrary to the choice of $H$.

Hereafter we may assume that $S$ is incident to fewer than $e(H)/2$ edges of $H$. Thus the minimum degree of a vertex in $S$ is at most $e(H)/\abs{S}=e(H)/(\gamma v(H))$. Removing edges incident to $S$ from $H$ then leaves a graph $H'$ with maximum degree at most $e(H)/\gamma v(H)$. In addition, since $e(H')\geq e(H)/2$ and $v(H')\leq v(H)$, the average degree of $H'$ is at least $e(H)/v(H)$.

By removing vertices of degree less than $e(H)/(4v(H))$ from $H'$, we obtain a graph $G_0$ on at least $v(H)/2 \ge v(H)^{1-\alpha/2}/2 \ge cn^{\alpha/2}$ vertices. Moreover, each vertex in this graph $G_0$ has degree between $e(H)/(4v(H))$ and $e(H)/(\gamma v(H))$.

\medskip\noindent\textit{Thinning the subgraph:} Note that the degrees of vertices in $G_0$ are about $e(H) / v(H)$, which is at least $cv(H)^\alpha \ge cm^\alpha$, where $m \eqdef v(G_0)$. We have already shown that $m \ge c n^{\alpha/2}$. However, it might happen that $e(H) / v(H)$ is much larger than the desired $cm^\alpha$. To rectify that, we shall find a decreasing sequence $d_0 \ge d_1 \ge \dots \ge d_k$ and a descending chain of subgraphs $G_0 \supseteq G_1 \supseteq \dots \supseteq G_k$ all on $m$ vertices, where $k = \lfloor\log_{2/3}\alpha\rfloor + 1$, such that
\begin{enumerate}
	\item $d_k = cm^\alpha$ and $d_i \ge (d_{i-1})^{2/3}$ for all $i \in \{1,2,\dots, k\}$;
	\item each vertex in the graph $G_i$ has degree between $\frac{1}{4}\left(\frac{2}{3}\right)^id_i$ and $\frac{1}{\gamma}\left(\frac{3}{2}\right)^id_i$ for all $i \in \{0,1,\dots,k\}$.
\end{enumerate}

There are many ways to define a desired decreasing sequence.
For example, we may first set $d_0 \eqdef e(H)/v(H)$ and then define recursively $d_i \eqdef \max((d_{i-1})^{2/3}, cm^\alpha)$. Since $d_0 < m$, our choice of $k$ ensures that $d_k = cm^\alpha$.

We obtain $G_{i+1}$ from $G_{i}$ by sampling each edge of $G_i$ with probability $p \eqdef d_{i+1} / d_{i}$. For each vertex $v$ of degree $d(v)$ in $G_{i}$, let $B(v)$ be the event that the degree of $v$ deviates from the expectation $pd(v)$ by at least $\frac{1}{3}pd(v)$ in the resulting random subgraph of $G_i$. Hoeffding's inequality shows that the probability $\Pr\left(B(v)\right)$ is at most $$2\exp\left(-2\left(\frac{p}{3}\right)^2d(v)\right) \le 2\exp\left(-\frac{2}{9}\left(\frac{d_{i+1}}{d_{i}}\right)^2 \frac{1}{4}\left(\frac{2}{3}\right)^id_i\right) = 2\exp\left(-\frac{1}{18}\left(\frac{2}{3}\right)^i\frac{(d_{i+1})^2}{d_{i}}\right).$$ As $d_{i+1} \ge (d_i)^{2/3}$, and $i\leq k-1=\lfloor \log_{2/3} \alpha\rfloor$, this means that
\[
  \Pr\left(B(v)\right)\leq 2\exp\left(-\frac{\alpha}{18} \left(d_i\right)^{1/3}\right).
\]
Since $d_i\geq cm^{\alpha}\geq c(cn^{\alpha/2})^\alpha$, and $n$ is sufficiently large, it follows that $$e \max \left\{d(v) : v \in V(G_i)\right\}\max \left\{\Pr\left(B(v)\right) : v \in V(G_i)\right\} < 1.$$
The symmetric version of the Lov\'asz local lemma then yields a subgraph $G_{i+1}$ of $G_i$ such that the degree of each vertex $v$ in $G_{i+1}$ is between $\frac{2}{3}pd(v)$ and $\frac{4}{3}pd(v) \le \frac{3}{2}pd(v)$. At the end, we obtain a graph $G' \eqdef G_k$ whose vertices have degrees between $\frac{\alpha}{6}cm^\alpha$ and $\frac{3}{2\alpha\gamma}cm^\alpha$.
\end{proof}
\end{document}